\documentclass[11pt,bezier]{article}
\usepackage{amsmath,graphicx,amssymb,amsfonts}

\newtheorem{prethm}{{\bf Theorem}}

\newenvironment{thm}{\begin{prethm}{\hspace{-0.5
               em}{\bf.}}}{\end{prethm}}
\newtheorem{prepro}{{\bf Theorem}}

\newenvironment{pro}{\begin{prepro}{\hspace{-0.5
               em}{\bf.}}}{\end{prepro}}
\newtheorem{precor}{{\bf Corollary}}

\newenvironment{cor}{\begin{precor}{\hspace{-0.5
               em}{\bf.}}}{\end{precor}}
\newtheorem{preconj}{{\bf Conjecture}}

\newtheorem{preremark}{{\bf Remark}}

\newenvironment{remark}{\begin{preremark}\em{\hspace{-0.5
               em}{\bf.}}}{\end{preremark}}
\newtheorem{prelem}{{\bf Lemma}}

\newenvironment{lem}{\begin{prelem}{\hspace{-0.5
               em}{\bf.}}}{\end{prelem}}
\newtheorem{preproof}{{\bf Proof.}}

\newenvironment{proof}[1]{\begin{preproof}{\rm
               #1}\hfill{$\Box$}}{\end{preproof}}

\title{\bf\large Some Relations between Rank, Chromatic\\ Number and
Energy of Graphs
}

\author{\small
{\bf \sc  S. Akbari}, {\bf \sc E. Ghorbani}, {\bf \sc S. Zare }
\thanks {AMS (2000) { \it Subject classification}: 05C15, 05C50,
15A03.}}
\date {}

\begin{document}
 \maketitle

\begin{abstract}
The energy of a graph $G$, denoted by $E(G)$, is defined as the
sum of the absolute values of all eigenvalues of  $G$. Let $G$ be
a graph of order $n$ and  ${\rm rank}(G)$ be the rank of the
adjacency matrix of $G$. In this paper we characterize all graphs
with $E(G)={\rm rank}(G)$. Among other results we show that apart
from a few families of graphs, $E(G)\geq 2\max(\chi(G),
n-\chi(\overline{G}))$, where $n$ is the number of vertices of
$G$, $\overline{G}$ and $\chi(G)$ are the complement and the
chromatic number of $G$, respectively. Moreover some new lower
bounds for $E(G)$ in terms of ${\rm rank}(G)$ are given.
 \end{abstract}
\vspace{1cm} \noindent{\bf\large Introduction} \vspace{4mm}

Let $G$ be  a graph. Throughout this paper the {\it order} of $G$
is the number of vertices of $G$. All the graphs that we consider
in the paper  are finite, undirected and simple. If $\{ v_1,
\ldots, v_n\}$ is the set of vertices of $G$, then the {\it
adjacency matrix} of $G$, $A=[a_{ij}]$,
 is an $n\times{n}$ matrix, where $a_{ij}=1$ if $v_i$ and $v_j$ are
 adjacent and $a_{ij}=0$
 otherwise. Thus $A$ is a symmetric matrix with zeros on the
 diagonal, and all eigenvalues of $A$ are real.  For a graph $G$,
   let ${\rm rank}(G)$ denote the rank
of the adjacency matrix of $G$. The {\it spectrum} of graph $G$,
 ${\rm Spec}(G)$, is the set of the eigenvalues of $A$, denoted by
 $\lambda_1(G)\geq \lambda_2(G)\geq \cdots \geq
 \lambda_n(G)$. We denote the {\it path} and the {\it complete graph}
of order $n$ by  $P_n$
  and $K_n$, respectively. The {\it
complete $t$-partite graph} is a graph whose vertices can be
partitioned into $t$ parts so that two vertices are adjacent if
and only if they belong to different subsets of the partition. We
 denote the complete $t$-partite graph with parts of sizes $r_1,
 \ldots, r_t$
  by $K_{r_1, \ldots, r_t}$. A {\it matching} of
$G$ is a set of  mutually non-incident edges. A {\it perfect
matching} of $G$ is a matching which covers all vertices of $G$.
For a graph $G$, {\it the chromatic number} of $G$, $\chi(G)$, is
the minimum number of colors needed to color the vertices of $G$
such that no two adjacent vertices have the same color.

 The {\it H\"{u}ckel \,molecular \,orbital}, HMO theory, is nowadays
 one
 of the
most important fields of theoretical chemistry where graph
eigenvalues occur. HMO theory deals with unsaturated conjugated
molecules. The vertices of the graph associated with a given
molecule are in one
 to one correspondence with the carbon atoms of the hydrocarbon system.
H\"{u}ckel theory in quantum chemistry insures that the total
$\pi$-electron energy of a conjugated hydrocarbon is simply the
energy of the corresponding molecular graph


 The {\it energy} of a graph $G$ is defined as the sum of the absolute
values of all
 eigenvalues and denoted by $E(G)$.
 If $\lambda_1, \ldots, \lambda_s$ are all the positive eigenvalues of
 a graph $G$, then we have $E(G)=2(\lambda_1+\cdots +
 \lambda_s)=-2(\lambda_{s+1}+\cdots + \lambda_n)$.   Recently much work on graph
energy appeared  in \cite{gut, rad, rati, ya1, ya2}.



\vspace{1cm} \noindent{\bf \large Main Results} \vspace{4mm}

\noindent First we state  the following lemma without proof.

\begin{lem}\label{k12} {\em\cite[p. 21]{spec}}
If for every eigenvalue $\lambda$ of a graph $G$, $\lambda\geq
-1$, then $G$ is a union of complete graphs.
\end{lem}

In \cite{faj}, it is shown that for any graph $G$, $E(G)\geq {\rm
rank}(G)$. Here we characterize all graphs $G$ for which $E(G)=
{\rm rank}(G).$

\begin{lem}\label{aa}  Let $G$ be a graph of order $n$. Then $E(G)\geq
{\rm rank}(G)$ and  equality holds if and only if
$G=\frac{r}{2}K_2\cup(n-r)K_1$ for some even positive integer $r$.
\end{lem}

\begin{proof}{
Assume that $\lambda_1,\ldots, \lambda_r$ are all non-zero
eigenvalues of $G$. Let
$\lambda^{n-r}(\lambda^r+a_1\lambda^{r-1}+\cdots+a_r)$  be the
characteristic polynomial of $G$, where  $a_r$ is a nonzero
integer. Then the arithmetic-geometric inequality implies that
\begin{equation}\label{arith}{\frac{|\lambda_1|+\cdots+|\lambda_r|}{r}\geq
\sqrt[r]{|\lambda_1|\cdots|\lambda_r|}=\sqrt[r]{|a_r|}\geq
1.}\end{equation}

 Thus $E(G)\geq {\rm rank}(G) $.
 If $G= \frac{r}{2}K_2 \cup(n-r)K_1$, obviously $E(G)={\rm rank}(G)$.
Conversely, suppose that $E(G)={\rm rank}(G)$. So equality holds
in (\ref{arith}), that is $|\lambda_1|=\cdots=|\lambda_r|=1$. Now,
by Lemma \ref{k12}, $G=\frac{r}{2}K_2\cup(n-r)K_1$.}
\end{proof}

  In \cite{bala}, it is shown that the energy of a connected graph with
 at least two vertices
   is greater than 1. In the following we improve this lower bound.
\begin{thm}\label{b}
For any connected graph $G$ apart from $K_1$ and $K_{1,i}, 1\leq
i\leq3$,
 $E(G)\geq4$.
\end{thm}

\begin{proof}{ We may assume that $G$ has at least four vertices.
 Clearly, ${\rm rank}(G)\geq2$.  If ${\rm rank }(G)=2$, it is
shown in \cite{cam} that $G$ is a complete bipartite graph. So $G$
is $K_{r,s}$, $rs\geq4$, and $E(G)\geq 2\sqrt{rs}\geq4$. If ${\rm
rank}(G)=3$, then $G$ has three non-zero eigenvalues.  Also $G$
has exactly one positive eigenvalue otherwise the sum of the
eigenvalues of $G$ is not zero because $\lambda_1(G)$ has the
greatest absolute value. Hence by Theorem 6.7 of \cite[p.
163]{spec}, $G$ is a complete multipartite graph. Since the rank
of a complete $t$-partite graph is $t$, $G$ is a complete
$3$-partite graph.
 Therefore $G$ has $K_3$
as an induced subgraph. Thus by Interlacing Theorem (Theorem
$0.10$ of \cite{spec}) $E(G)\geq E(K_3)=4$. If $ {\rm rank
}(G)\geq 4$, then by Lemma \ref{aa}, $E(G)\geq4$.}
\end{proof}

\begin{thm}\label{gho}
If $G$ is a connected bipartite graph of rank $r$, then
$E(G)\geq\sqrt{(r+1)^2-5}$.
\end{thm}
\begin{proof}{ Let $\lambda_1,\ldots,\lambda_s$ be the positive
eigenvalues of $G$, where $s=r/2$. Then
$$E^2(G)=\left(2\sum_{i=1}^s{\lambda_i}\right)^2=
4\left(\sum_{i=1}^s{\lambda_i}^2+\sum_{i\neq{j}}{\lambda_i}{\lambda_j}\right)=
4\left(m+s(s-1)a\right),$$ where $m$ is the number of edges in $G$
and $a$ is the arithmetic mean of
$\{\lambda_i\lambda_j\}_{i\neq{j}}$. The geometric mean of
$\{\lambda_i\lambda_j\}_{i\neq{j}}$ is

$$\left(\prod_{i\neq{j}}\lambda_i
\lambda_j\right)^{(s(s-1))^{-1}}=k^{1/s},$$ where
$k={\lambda_1}^2\cdots{\lambda_s}^2$. Since $G$ is connected,
$m\geq{r-1}$. Note that $k\geq1$. So we have
$$E(G)\geq\sqrt{4m+r(r-2)\sqrt[r]{k^2}}\geq{\sqrt{(r+1)^2-5}}.$$}

\end{proof}

\noindent The proof of the following lemma is easy and we leave it
to the reader.

\begin{lem}\label{akb}
If $T$ is a tree with no perfect matching and isolated vertex,
then $T$ has at least two maximum matchings.
\end{lem}

The following lemma is an immediate consequence of Harary's
Theorem, see \cite[p. 44]{big}.

\begin{lem}\label{har} The number of maximum matchings of a tree is
equal to the product of its non-zero eigenvalues.
\end{lem}

\begin{cor}\label{2} If $T$ is a tree with no perfect matching,
 then the product of its non-zero eigenvalues is at least $2$.
 \end{cor}

\begin{thm}\label{d}
Let $G$ be a  bipartite graph with at least $4$ vertices. If $G$
is not  full rank, then $E(G)\geq 1+{\rm rank}(G)$.
\end{thm}

\begin{proof}{Without loss of generality, we may assume that $G$ is a
 connected
graph. From the proof of Theorem \ref{gho}, we have
$E(G)\ge\sqrt{4m+r(r-2)\sqrt[r]{k^2}}$, where $r={\rm rank}(G)$
and $k=\lambda_1^2\cdots\lambda_s^2$ and $\lambda_1, \ldots,
\lambda_s$ are positive eigenvalues of $G$, $s=r/2$. If $G$ is a
tree, then by Theorem 8.1 of  \cite{spec}, $G$ has no perfect
matching. Thus Corollary \ref{2} implies that $k\geq 2$. Hence
$E(G)\geq \sqrt{4(n-1)+ r(r-2)\sqrt[r]4}$. Note that $G$ is not
full rank. So we have $E(G)\geq \sqrt{4r+ r(r-2)\sqrt[r]4}$.
 Note that if $r\geq3$, then
 $\sqrt[r]4>\exp(1/r)>1+\frac{1}{r}\ge1+\frac{1}{r(r-2)}$.
 Thus if $r\geq3$, $E(G)>\sqrt{4r+r(r-2)+1}=r+1$.
 If $r=2$, then by Theorem \ref{b} we are done. If
$G$ is not a tree, then since $G$ is not  full rank, we find,
$m\ge n\ge r+1$ and the proof is complete.}\end{proof}


Now, we would like to obtain some lower bounds for  $E(G)$  in
terms of the chromatic number of $G$ and the chromatic number of
$\overline{G}$.

\begin{pro}\label{a} {\rm (Theorem 2.30 of \cite{fav})}
For any graph $G$,  $n-\chi(\overline{G})\leq
\lambda_1+\cdots+\lambda_{\chi(\overline{G})}.$
\end{pro}

\noindent By Theorem \ref{a}, we have the following result.

\begin{thm} \label{main}  For every graph $G$,
$E(G)\geq 2 (n-\chi(\overline{G})).$
\end{thm}

\begin{remark} \label{rem} A well-known theorem of Nordhaus and Gaddum
\cite{nor} states that for every graph $G$ of order $n$,
$\chi(G)+\chi(\overline{G})\leq n+1$. The graphs attaining
equality in the Nordhaus-Gaddum Theorem were characterized by
Finck \cite{fin}, who proved that there are exactly two types of
such graphs, the types $(a)$ and $(b)$ defined as follows.

 \noindent (i) A graph
$G$ is of type $(a)$ if it has a vertex $v$ such that
$V\setminus\{ v\}$ can be partitioned into subsets $K$ and $S$
with the properties that $K\cup \{v\}$ induces a clique of $G$ and
$S\cup \{v\}$ induces an independent set of $G$ (adjacency between
$K$ and $S$ is arbitrary). Note that if $G$ is of type $(a)$, then
so does its complementary graph $\overline{G}$.

\noindent  (ii) A graph $G$ is of type $(b)$ if it has a subset
$C$ of five vertices such that $V\setminus C$ can be partitioned
into subsets $K$ and $S$ with the properties that $K$ induces a
clique, $S$ induces an independent set, $C$ induces a $5$-cycle,
and every vertex of $C$ is adjacent to every vertex of $K$ and to
no vertex of $S$(adjacency between $K$ and $S$ is arbitrary). Note
that if $G$ is of type $(b)$, then so does its complementary
graph.
\end{remark}

If we omit a perfect matching from the complete graph $K_{2n}$,
the resulting graph is called {\it cocktail party} and denoted by
$CP(n)$. For any graph $G$ with vertices $\{v_1,\ldots,v_n\}$, and
any non-negative integers $a_1,\ldots,a_n$, we construct the {\it
generalized line graph} $L(G;a_1,\ldots,a_n)$ as follows:

 The
vertex set is the union of the vertex sets of $L(G)$, $CP(a_1),
\ldots, CP(a_n)$, and the edge set is the union of the edge sets,
together with edges joining all vertices of $CP(a_i)$ to every
vertex of $L(G)$ corresponding to an edge of $G$ containing $v_i$,
for $1\le i\le n$.


Denote by $A_{n,t}$ for $1\leq t\leq n-1$ the graph obtained by
joining a new vertex to $t$ vertices of the complete graph $K_n$.
If we add two pendant vertices to a common vertex of $K_n$, then
the resulting graph has order $n+2$ and we denote it  by $B_n$.
For the proof of the next theorem we need the following
interesting result due to Wilf, see \cite[p. 55]{big}.

\begin{lem} \label{chrm} For any graph $G$, $\chi(G)\leq
\lambda_1(G)+1$, where $\lambda_1(G)$ denotes the largest
eigenvalue of $G$.
\end{lem}

\begin{thm} \label{ab} Let $G$ be a graph. Then $E(G)<2\chi(G)$ if and
only if $G$ is a union of some isolated vertices and one of the
following graphs:\\ (i) the complete graph $K_n$;\\ (ii) the graph
$B_n$;\\
(iii) the graph $A_{n,t}$ for $n\leq7$, except when $(n,t)=(7,4)$,
and for $n\geq8$ with $t\in \{1, 2, n-1 \}$; \\ (iv) a triangle
with two pendant vertices adjacent to different vertices (see the
graph $H_5$ in Figure {\rm \ref{345}}).
\end{thm}
\begin{proof}{ First we show that the graphs stated in the theorem
satisfy $E(G)<2\chi(G)$. The assertion is clear for complete
graphs. The characteristic polynomial of $B_n$ (see \cite[p.
159]{spec}) is
 $$\lambda(\lambda+1)^{n-2}[(\lambda^2-2)(\lambda-n+2)
 +\lambda(\lambda-n+1)(\lambda+1)+\lambda^2(\lambda-n+2)].$$
 Therefore  $B_n$ has at least $n-2$
 eigenvalues $-1$.
   By Lemma \ref{k12} we find $\lambda_{n+2}<-1$, so by
 Theorem $6.7$ of \cite{spec},
  $B_n$  has exactly two positive eigenvalues.
This implies that
 $E(B_n)=-2(\lambda_{n+2}-n+2)$ since the sum of eigenvalues of any
 graph is zero.
  On the other hand the characteristic polynomial of $B_n$ is
 $\lambda(\lambda+1)^{n-2}f(\lambda)$, where
$f(\lambda)=\lambda^3+(2-n)\lambda^2-(1+n)\lambda+2n-4$.
  It is not hard to see that  $f(\lambda)<f(-2)=-2$ for any
 $\lambda<-2$. Therefore all
eigenvalues of $B_n$ are  more
  than $-2$. Thus $E(B_n)=-2(\lambda_{n+2}-n+2)<2n=2\chi(B_n)$.

 A calculation shows  that for $n\leq7$, $E(A_{n,t})<2n=2\chi(A_{n,t})$
 except
$E(A_{7,4})=14$. So we may assume that $n\geq 8$. The graph $K_n$
has $n-1$ eigenvalues $-1$, therefore by Interlacing Theorem, the
graph $A_{n,t}$ has at least $n-2$ eigenvalues $-1$. On the other
hand the graphs $A_{n,1}$ and $A_{n,2}$ are not complete
multipartite graphs, so they have at least two positive
eigenvalues. Then  again Interlacing Theorem implies that these
two graphs have exactly two positive eigenvalues. The graphs
$A_{n, 1}$ and $A_{n, 2}$ are line graphs, and $A_{n, n-1}=L(K_{1,
n-1};1,0,\ldots,0)$, where the vertex $v_1$ in $K_{1,n}$ is the
vertex with maximum degree. Hence their eigenvalues are at least
$-2$ and Theorem 1.6 of \cite{recent} shows that
$\lambda_{n+1}>-2$. Thus $E(A_{n,t})=-2(\lambda_{n+1}-n+2)<2n$ for
$t=1,2$. The graph $A_{n,n-1}$ has a zero eigenvalue, and in the
same way we find that $E(A_{n,n-1})<2n$. Now, let $2<t<n-1$. We
claim that the graph $A_{n,t}$ cannot be a generalized line graph.
For any $l\geq 2$, $CP(l)$ has $C_4$ as an induced subgraph and if
$A_{n,t}$ is a generalized line graph, then we conclude that
$A_{n,t}=L(K_{1,n-1};1)$. Clearly, $A_{n,t}$ is not isomorphic to
$L(K_{1,n-1};1)$. Thus by Exercise 14 of \cite[p. 278]{god},
$\lambda_{n+1}\leq-2$. Hence $E(A_{n,t}) \geq 2n$. Finally for the
graph $H_5$ (see Figure \ref{345}) by an easy computation we see
$E(H_5)<2\chi(H_5)$ (see \cite{rea}).

 Now, we show that apart from exceptional  cases of the theorem, for
 any graph $G$,
 $E(G)\geq 2\chi(G)$. If $G$
has two non-trivial components, then $G$ has $2K_2$ as an induced
subgraph. Hence Interlacing Theorem (Theorem 0.10 of \cite{spec})
and Lemma \ref {chrm} imply that  $E(G)\geq 2(\lambda_1+
\lambda_2) \geq 2(\lambda_1+1) \geq 2\chi(G)$.
 Since  isolated vertices do not contribute to the energy and the
chromatic
 number, we may assume that $G$ is connected. By Remark \ref{rem}, we
 have
$\chi(G)+\chi(\overline{G})\leq n+1$. If
$\chi(G)+\chi(\overline{G})\leq n$, then we are done by Theorem
\ref{main}. So we let $\chi(G)+\chi(\overline{G})=n+1$. In this
case $G$ is either of type $(a)$ or of type $(b)$. If $G$ is of
type $(b)$, then $G$ has $C_5$ as an induced subgraph. Therefore
$\lambda_2(G)+\lambda_3(G)\geq\lambda_2(C_5)+\lambda_3(C_5)>1$.
Note that $\lambda_3(C_5)>0$.  Thus by Lemma \ref{chrm} we have
$$E(G)\geq 2(\lambda_1+\lambda_2+\lambda_3)>2(1+\lambda_1)\geq
2\chi(G).$$

Thus one may assume that $G$ is of type $(a)$. For simplification
let $|K|=t$, where $K$ is a complete subgraph of $G$ defined in
Remark \ref{rem}. It is easily seen that $\chi(G)=t+1$. Clearly,
$K_{t+1}$ is an induced subgraph of $G$. We know that $K_{t+1}$
has one eigenvalue $t$ and $t$ eigenvalues $-1$. So by Interlacing
Theorem,
 $G$ has at least $t$ eigenvalues which are at most $-1$. If $G$ has an
induced subgraph with at least one eigenvalue $\lambda$ such that
$\lambda \leq -2$, then the sum of all negative eigenvalues of $G$
is less than $-(t-1)-2=-t-1$. Thus $E(G)\geq 2t+2=2\chi(G).$

Therefore we may assume that every eigenvalue of each induced
subgraph of $G$ is more than $-2$. This implies that $G$ has no
$K_{1,4}$ as an induced subgraph. Hence every vertex of $K$ is
adjacent to at most two vertices of $S$.

 First suppose that there is
a vertex $a \in K$ which is adjacent to two vertices $\{ x, y\}
\subseteq S$. If $|S| \geq 3$, then there exists a vertex $z \in
S\backslash \{x , y\}$ such that $z$ is adjacent to a vertex $b
\in K$ and $b\neq a, v$, where $v$ is the vertex given in Remark
\ref{rem}. Thus $G$ has either $H_1$, or $H_2$ as an induced
subgraph (see Figure \ref{12}).

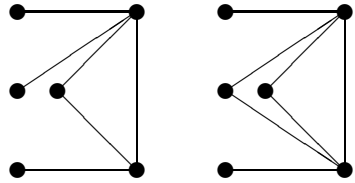
\begin{figure}[h]\begin{center}
\setlength{\unitlength}{3pt}
\begin{picture}(25,25)(0,0)
\put(5,20){\circle*{2}}
 \put(20,20){\circle*{2}}
\put(5,0){\circle*{2}} \put(20,0){\circle*{2}}
\put(10,10){\circle*{2}} \put(5,10){\circle*{2}}
\put(5,20){\line(1,0){15}} \put(20,20){\line(-1,-1){10}}
\put(20,20){\line(0,-1){20}} \put(5,0){\line(1,0){15}}
\put(20,0){\line(-1,1){10}}\put(5,10){\line(3,2){15}}
\end{picture}
\begin{picture}(25,25)(0,0)
\put(5,20){\circle*{2}}
 \put(20,20){\circle*{2}}
\put(5,0){\circle*{2}} \put(20,0){\circle*{2}}
\put(10,10){\circle*{2}} \put(5,10){\circle*{2}}
\put(5,20){\line(1,0){15}} \put(20,20){\line(-1,-1){10}}
\put(20,20){\line(0,-1){20}} \put(5,0){\line(1,0){15}}
\put(20,0){\line(-1,1){10}}\put(5,10){\line(3,2){15}}
\put(5,10){\line(3,-2){15}}
\end{picture}\end{center}\caption{The graphs $H_1$ and $H_2$}\label{12}
\end{figure}



\noindent We have (see \cite{rea}) $\lambda_6(H_1)<-1.8$,
$\lambda_5(H_1)<-1.3$, and $\lambda_6(H_2)<-1.7$,
$\lambda_5(H_2)<-1.6$. Since $G$ has at least $t$ eigenvalues
which are at most $-1$, the sum of all negative eigenvalues of $G$
is less than $-(t-2)-3=-t-1$. Thus $E(G)\geq 2t+2=2\chi(G).$ If
$|S|=2$, there exist two cases: (1) $G=B_{n-2}$; (2) $G$ has
either $K_{1,1,3}$ or $H$ as an induced subgraph, where $H$ is the
graph obtained by removing one of the pendant vertices of $H_2$.
If (1) is the case, then $E(G)<2\chi(G)$. If (2) is the case,
since $\lambda_5(K_{1,1,3})=-2$, $\lambda_5(H)<-1.74$, and
$\lambda_4(H)<-1.27$, as before we conclude that
$E(G)\ge2\chi(G)$.

 Now, suppose that every
vertex in $K$ is adjacent to at most one vertex of $S$. If $|S|
\geq2$ and $ | K|\geq 3$, then $G$ has an induced subgraph
isomorphic to $H_3$ or $H_4$.\\

\begin{figure}[h]\begin{center}
\setlength{\unitlength}{2.5pt}
\begin{picture}(25,25)(0,0)
\put(0,25){\circle*{2}}\put(10,25){\circle*{2}}\put(25,25){\circle*{2}}
\put(0,0){\circle*{2}}\put(10,0){\circle*{2}}\put(25,0){\circle*{2}}
\put(10,25){\line(1,0){15}}\put(10,25){\line(0,-1){25}}\put(10,25){\line(3,-5){15}}
\put(10,25){\line(-1,0){10}}
\put(10,0){\line(1,0){15}}\put(10,0){\line(-1,0){10}}\put(10,0){\line(3,5){15}}
\put(25,0){\line(0,1){25}}
\end{picture}\hspace{1cm}
\begin{picture}(25,25)(0,0)
\put(2.5,12.5){\circle*{2}}\put(10,25){\circle*{2}}\put(25,25){\circle*{2}}
\put(10,0){\circle*{2}}\put(25,0){\circle*{2}}\put(32.5,12.5){\circle*{2}}
\put(10,25){\line(1,0){15}} \put(10,25){\line(0,-1){25}}
\put(10,25){\line(3,-5){15}} \put(2.5,12.5){\line(3,5){7.5}}
\put(10,0){\line(1,0){15}} \put(10,0){\line(3,5){15}}
\put(25,0){\line(0,1){25}} \put(32.5,12.5){\line(-3,5){7.5}}
 \put(32.5,12.5){\line(-3,-5){7.5}}
\end{picture}
\hspace{1cm}
\begin{picture}(25,25)(0,0)
\put(5,25){\circle*{2}}
 \put(25,25){\circle*{2}}
\put(5,0){\circle*{2}} \put(25,0){\circle*{2}}
\put(12.5,12.5){\circle*{2}} \put(5,25){\line(1,0){20}}
\put(25,25){\line(-1,-1){12.5}} \put(25,25){\line(0,-1){25}}
\put(5,0){\line(1,0){20}} \put(25,0){\line(-1,1){12.5}}
\end{picture}
\end{center}\caption{The graphs $H_3$, $H_4$, and $H_5$}\label{345}
\end{figure}
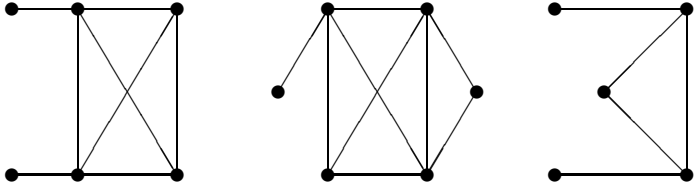

 The graph $H_3$ has two eigenvalues one of which is less
than $-1.39$ and  the other one is less than $-1.61$. Also $H_4$
has two eigenvalues one of which is less than $-1.3$ and  the
other one is less than $-1.7$, see \cite{rea}. Thus $E(G)\geq 2
\chi(G)$. It remains to consider the case $|S|\leq 1$ or $|K|\leq
2$. If $|S|= 0$, then $G=K_n$. If $|S|= 1$, then $G= A_{n,t}$ for
some $t$. If $|K| \leq 2$, then $|S|\leq 2$. It can be easily
checked that $G$ is one of the graphs $B_1=K_{1,2}$, $A_{3,1}$,
$A_{3,2}$, $H_5$ or $K_i$, $i=1,2,3$.}\end{proof}


Let $G$ be a connected graph of order $n$. The following corollary
shows that either the graph $G$, or $\overline{G}$ has energy at
least $n$. Compare with Corollary 5.2 of \cite{gut} which states
that if $G$ has no zero eigenvalue, then $E(G)\geq n$.

\begin{cor} Let $G$ be a graph of order $n\geq3$. If  $G$ or
 $\overline{G}$
is neither a complete graph nor one of the graphs $A_{k,k-1}$,
 $B_1$, $B_2$, and $A_{3,1}$, then
 $E(G)+ E(\overline{G})\geq2n.$
\end{cor}
\begin{proof}{ If $\overline{G}$ is not one of the graphs
described in Theorem \ref{ab}, then $E(\overline{G})\ge
2\chi(\overline{G})$ and the corollary follows from Theorem
\ref{main}. If $G=H_5$, then $E(G)+E(\overline{G})>10$. If $G$ is
a complete graph, $B_1$, $B_2$, or $A_{3,1}$, it is easily seen
that $E(G)+ E(\overline{G})< 2n$. If $G= A_{k,k-1}$, then
$\overline{G}= K_2\cup (k-1)K_1$. Hence
$E(A_{k,k-1})+E(\overline{A_{k,k-1}})<2k+2=2n$. To complete the
proof, it is enough to show that the theorem holds for $B_k$,
$k\geq3$ and $A_{k,t}$ for $k\geq4$, and $t=1,2$. The graph $B_k$
has $k-2$ eigenvalues $-1$, and
$\lambda_n(B_k)\leq\lambda_6(B_4)<-1.8$. So $E(B_k)> 2(k-0.2)$. On
the other hand $K_{1,1,2}$ is an induced subgraph of
$\overline{B_k}$. Therefore $E(\overline{B_k})\geq E(K_{1,1,2})>5$
and so $E(B_k)+E(\overline{B_k})\geq 2(k+2).$

 The graph $A_{k,t}$ has $k-1$ eigenvalues $-1$ and
$\lambda_n(A_{k,1})\leq\lambda_5(A_{4,1})<-1.5$, so
$E(A_{k,1})>2k-1$. Also
$\lambda_n(A_{k,2})\leq\lambda_5(A_{4,2})<-1.68$, hence
$E(A_{k,2})>2k-0.64$. Now, the facts $E(\overline{A_{k,1}})\geq
E(K_{1,3})>3.4$ and $E(\overline{A_{k,2}})\geq E(K_{1,2})>2.8$
complete the proof.}

\end{proof}

{\bf Acknowledgements. }The first and second authors are indebted
to the Institute for Studies in Theoretical Physics and
Mathematics (IPM) for support. Also the first and third authors
are indebted to the Abdus Salam  International Center for
Theoretical Physics (ICTP) for hospitality during their visits.
The research of the first author was in part supported by a grant
(No.\,85160211) from IPM. The authors would like to thank the
referees for their fruitful comments and suggestions.

\noindent {\small {\sc S. Akbari\,${}^{a, b}$    \quad\quad \quad
 \quad\,\,\, {\tt s\_akbari@sharif.edu}}\\
\noindent{\sc E. Ghorbani\,${}^b$  \quad \quad   \quad \quad\,{\tt
e\_ghorbani@math.sharif.edu
}} \\
\noindent{\sc S. Zare\,${}^b$  \quad  \quad \quad  \quad \quad
\quad\, {\tt
sa\_zare\_f@yahoo.com}} \\
\noindent ${}^a$\,Institute for Studies in Theoretical Physics and
Mathematics, Tehran, Iran. \quad \\
${}^b$\,Department of Mathematical Sciences, Sharif University of
Technology,\\\noindent P. O. Box 11365-9415, Tehran, Iran.}
\end{document}